\documentclass[10pt]{article}
\usepackage{amsfonts}

\usepackage{mathrsfs}
\usepackage{amsmath,url}
\usepackage{mathrsfs,amscd,amssymb,amsthm,amsmath,bm,graphicx,psfrag,subfigure}

\makeatletter

\renewcommand{\@seccntformat}[1]{{\csname the#1\endcsname}{\normalsize .}\hspace{.5em}}
\makeatother

\def \[{\begin{equation}}
\def \]{\end{equation}}

\newtheorem{thm}{Theorem}[section]

\newtheorem{lem}[thm]{Lemma}

\newenvironment{wst}
{\setlength{\leftmargini}{1.5\parindent}
 \begin{itemize}
 \setlength{\itemsep}{-1.1mm}}
{\end{itemize}}

\begin{document}
\setlength{\baselineskip}{15pt}
\begin{center}{\Large \bf On the eccentric distance sum of unicyclic graphs with a given\\[3pt] matching number\footnote{Financially supported by the National Natural Science Foundation of China (Grant Nos. 11071096, 11271149).}}

\vspace{4mm}

{\large Shuchao Li$^{a,}$\footnote{Corresponding author.\ E-mail: lscmath@mail.ccnu.edu.cn (S.C.
Li), 425333559@qq.com (Y.B. Song), bwei@olemiss.edu (B. Wei)},\ \ Yibing Song$^a$,\ \ Bing Wei$^b$}\vspace{2mm}

$^a$Faculty of Mathematics and Statistics,  Central China Normal
University, Wuhan 430079, P.R. China\vspace{2mm}

$^b$Department of Mathematics, University of Mississippi University, Oxford, MS 38677, USA
\end{center}
\noindent {\bf Abstract}: Let $G = (V_G,E_G)$ be a simple connected graph. The eccentric distance sum of $G$ is defined as $\xi^d(G)=\sum_{v \in V_G}\,\varepsilon_G(v)D_G(v),$ where $\varepsilon_G(v)$ is the eccentricity of the vertex $v$ and $D_G(v)=\sum_{u \in V_G}\,d(u,v)$ is the
sum of all distances from the vertex $v$. In this paper, we characterize $n$-vertex unicyclic graphs with given matching number having the minimal and second minimal eccentric distance sums, respectively.

\vspace{2mm} \noindent{\it Keywords}: Eccentric distance sum; Unicyclic graph; Matching number

\vspace{2mm}


 {\setcounter{section}{0}
\section{\normalsize Introduction}\setcounter{equation}{0}
We consider only simple connected graphs (i.e. finite, undirected graphs without loops or multiple edges). For a simple $G = (V_G, E_G)$ with $u,\,v \in V_G$, let $d_G(u)$ (or $d(u)$ for short) denote the degree of $u$ and $\Delta(G)=\max \{d(x)|x \in V_G\}$; $N_G(u)$ (or $N(u)$ for short) denotes the set of all the adjacent vertices of $u$ in $G$ and $N_G[u]=N_G(u)\cup \{u\};$ the \textit{distance} $d_G(u,v)$ is defined as the length of a shortest path connecting vertices $u$ and $v$ in $G$. The \textit{eccentricity} of $u,$ denoted by $\varepsilon_G(u)$, is the maximum distance from $u$ to all other vertices in $G.$

Two distinct edges in a graph $G$ are \textit{independent} if they do not have a common endvertex in $G.$ A set of pairwise independent edges of $G$ is called a
\textit{matching} of $G,$ while a matching of maximum cardinality is a \textit{maximum matching} of $G.$ Let $M$ be a matching of $G.$ The vertex $v$ in $G$ is $M$-\textit{saturated} if $v$ is incident with an edge in $M;$ otherwise, $v$ is $M$-\textit{unsaturated}. A \textit{perfect matching} $M$ of $G$ means that each vertex of $G$ is $M$-saturated; clearly, every perfect matching is maximum. The \textit{matching number} $m$ of $G$ is the cardinality of a maximum matching of $G;$ we also call such maximal matching as an $m$-\textit{matching.} An $M$-\textit{alternating path} of $G$ is a path whose edges are alternately in $E\backslash M$ and $M.$ An $M$-\textit{augmenting path} is an $M$-alternating path whose origin and terminus are $M$-unsaturated. We call $G$ a \textit{unicyclic graph,} if $G$ is connected satisfying $|V_G|=|E_G|.$

A single number that can be used to characterize some property
of the graph of a molecule is called a \textit{topological index},
or \textit{graph invariant}. Topological index is a graph theoretic
property that is preserved by isomorphism. The chemical information
derived through topological index has been found useful in chemical
documentation, isomer discrimination, structure property correlations,
etc. \cite{0}. For quite some time there has been rising interest
in the field of computational chemistry in topological indices. The
interest in topological indices is mainly related to their use in
nonempirical quantitative structure-property relationships and
quantitative structure-activity relationships.

Among distance-based topological indices,
the \textit{Wiener index} has been one of the most widely used descriptors in
quantitative structure activity relationships, which is defined as the sum of all distances between unordered pairs of vertices
$
  W(G) =\sum_{u,v\in V_G}d(u, v).
$
It is considered as one of the most used topological index with high correlation with many physical and
chemical properties of a molecule (modelled by a graph). For the recent survey on Wiener index one may refer to \cite{4}
and the references cited in.

Another distance-based topological index is the  \textit{degree distance index} $DD(G),$ which was introduced by Dobrynin and Kochetova \cite{05} and
Gutman \cite{09} as graph-theoretical descriptor for characterizing alkanes; it can be considered as a weighted version of the
Wiener index
$$
DD(G) = \sum_{u,v\in V_G}(d_G(u)+ d_G(v))d(u, v) = \sum_{v\in V_G}d_G(v)D_G(v),
$$
where the summation goes over all pairs of vertices in $G$.

Recently, a novel distance-based topological index is the \textit{eccentric distance sum} (EDS), which was introduced
by Gupta, Singh and Madan \cite{S0}. It is defined as
$$
\xi^d(G) =\sum_{v\in V_G} \varepsilon(v)D_G(v)=\sum_{u,v\in V_G} (\varepsilon(v)+\varepsilon(u))d(u,v).
$$
This topological index has vast potential in structure activity/property relationships; it also displays
high discriminating power with respect to both biological activity and physical properties; see \cite{S0}. From \cite{S0} we
also know  that some structure activity and quantitative structure property
studied using \textit{eccentric distance sum} were better than the corresponding
values obtained using the Wiener index. It is also interesting to study the
mathematical property of this topological index.
Yu, Feng and Ili\'{c} \cite{12} identified the
extremal unicyclic graphs of given girth having the minimal and second minimal EDS;
they also characterized the trees with the minimal EDS among the $n$-vertex trees of a given diameter.
Hua, Xu and Shu \cite{H-H-X-K}
obtained the sharp lower bound on EDS of $n$-vertex cacti. Hua, Zhang and Xu \cite{3} studied the relationship of
EDS with some other graph parameters over some typical class of graphs.
Ili\'{c}, Yu and Feng  \cite{13} studied the various lower and upper bounds for
the EDS in terms of the other graph invariant including the Wiener index,
the degree distance index and so on.
Zhang, Yu, Feng and one of the present authors \cite{17} determined the minimal and second minimal eccentricity distance sums among
all the $n$-vertex trees with mathching number $m$.

In this paper, we study the EDS of unicyclic graphs with given matching number. We determine the minimal and second minimal eccentricity distance sums among
all the $n$-vertex unicyclic graphs with mathching number $m$. The corresponding extremal graphs are also characterized.

\section{\normalsize  Preliminary Results }\setcounter{equation}{0}
Given a unicyclic graph $G$ with a unique cycle $C_k=u_1u_2\dots u_ku_1,$ $G-E_{C_k}:=T_1\cup T_2 \cup \cdots \cup T_k,$ where $T_i$ is a tree containing the vertex $u_i$ as its \textit{root}; and $T_i$ is called the \textit{pendant tree} of $G$, $i=1,2,\dots,k.$

For integers $n \geq 3$ and $1 \leq m \leq \lfloor\frac{n}{2}\rfloor$, let $\mathscr{U}_{n,m}$ be the set of $n$-vertex unicyclic graph with matching number $m.$ Obviously, $\mathscr{U}_{n,1}=\{C_3\}$. In what follows, we assume that $2 \leq m \leq \lfloor\frac{n}{2}\rfloor$. And denote $\mathscr{U}^2_{n,m}=\{G \in \mathscr{U}_{n,m}|\varepsilon_G(c)=2,\,c \,\,\hbox{is a central-vertex of $G$}\}$. Furthermore, we need the following notation:
\begin{itemize}
  \item Let $C_3(1^{k_1}, 1^{k_2}, 1^{k_3})$ be the graph obtained from $C_3=u_1u_2u_3u_1$ by inserting $k_i$ pendant vertices at $u_i,\, i=1,2,3$. Let $H(n,3;1^{k_1}S_{t_1,t_2,\dots,t_r},1^{k_2},1^{k_3})$ be the $n$-vertex graph obtained from $C_3(1^{k_1},1^{k_2},1^{k_3})$ and $r(r \geq 1)$ stars, say $K_{1,t_1+1},K_{1,t_2+1},\dots,K_{1,t_r+1},$ by identifying $u_1$ of $C_3(1^{k_1},1^{k_2},1^{k_3})$ with just one pendant vertex of each of the $r$ stars, respectively (see Fig. 1). In particular, if $t_1=t_2=\cdots=t_r=1,$ then denote the resultant graph by $H(n,3;1^{k_1}2^r,1^{k_2},1^{k_3}).$
  \item Let $C_4(1^{k_1},1^{k_2},1^{k_3})$ be the graph obtained from $C_4=u_1u_2u_3u_4u_1$ by inserting $k_i$ pendant vertices at $u_i,\, i=1,2,3$. Let $H(n,4;1^{k_1},1^{k_2}S_{t_1,t_2,\dots,t_r},1^{k_3})$ be the $n$-vertex graph obtained from $C_4(1^{k_1},1^{k_2},1^{k_3})$ and $r(r \geq 1)$ stars, say $K_{1,t_1+1},K_{1,t_2+1},\dots,K_{1,t_r+1},$ by identifying $u_2$ of $C_4(1^{k_1},1^{k_2},1^{k_3})$ with just one pendant vertex of each of the $r$ stars, respectively (see Fig. 1). In particular, if $t_1=t_2=\cdots=t_r=1,$ then denote the resultant graph by  $H(n,4;1^{k_1},1^{k_2}2^r,1^{k_3}).$
  \item Let $C_5(1^{k_1}, 1^{k_2}, 1^{k_3})$ be the graph obtained from $C_5=u_1u_2u_3u_4u_5u_1$ by inserting $k_i$ pendant vertices at $u_i,\, i=1,2,3$. Let $H(n,5;1^{k_1},1^{k_2}S_{t_1,t_2,\dots,t_r},1^{k_3})$ be the $n$-vertex graph obtained from $C_5(1^{k_1},1^{k_2},1^{k_3})$ and $r(r \geq 1)$ stars, say $K_{1,t_1+1},K_{1,t_2+1},\dots,K_{1,t_r+1},$ by identifying $u_2$ of $C_5(1^{k_1},1^{k_2},1^{k_3})$ with just one pendant vertex of each of the $r$ stars, respectively (see Fig. 1). In particular, if $t_1=t_2=\ldots=t_r=1,$ then denote the resultant graph by  $H(n,5;1^{k_1},1^{k_2}2^r,1^{k_3}).$
\end{itemize}
\begin{figure}[h!]
\begin{center}
  \psfrag{e}{$k_2$}\psfrag{A}{$H(n,5;1^{k_1},1^{k_2}S_{t_1,t_2,\dots,t_r},1^{k_3})$}
  \psfrag{f}{$k_3$}\psfrag{B}{$H(n,3;1^{k_1}S_{t_1,t_2,\dots,t_r},1^{k_2},1^{k_3})$}
  \psfrag{C}{$H(n,4;1^{k_1},1^{k_2}S_{t_1,t_2,\dots,t_r},1^{k_3})$}
  \psfrag{g}{$m-3$}\psfrag{a}{$t_1$}\psfrag{b}{$t_2$}
   \psfrag{p}{$u_1$}\psfrag{q}{$u_2$}\psfrag{r}{$u_3$}\psfrag{s}{$u_4$}\psfrag{t}{$u_5$}
  \psfrag{c}{$t_r$}\psfrag{d}{$k_1$}
  \includegraphics[width=150mm]{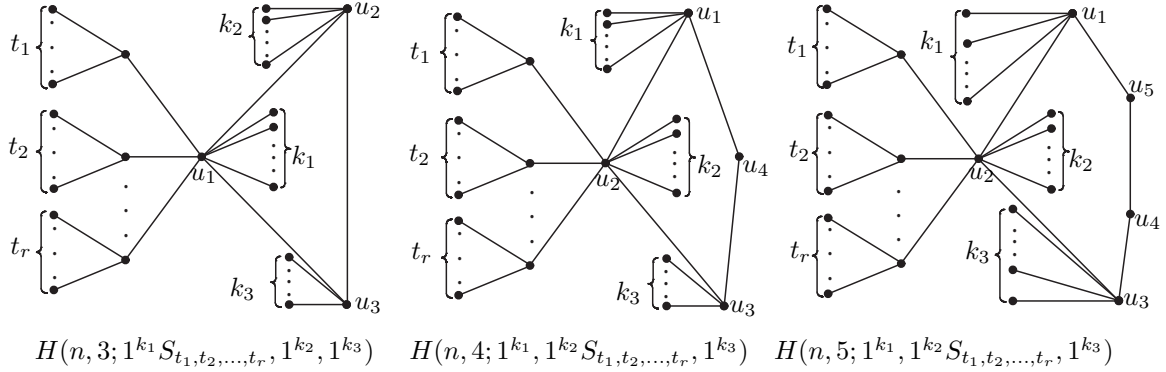}\\
  \caption{Graphs $H(n,3;1^{k_1}S_{t_1,t_2,\dots,t_r},1^{k_2},1^{k_3})$ and $H(n,i;1^{k_1},1^{k_2}S_{t_1,t_2,\dots,t_r},1^{k_3}),\, i=4,5.$}
\end{center}
\end{figure}

For convenience, let
\begin{equation*}\label{eq:2.1}
\begin{array}{ll}
  U_{n,m}:=H(n,3; 1^{n-2m+1}2^{m-2},1^0,1^0), &  U_{n,m}':=H(n,5;1^0,1^{n-2m+1}2^{m-3},1^0),\\[5pt]
   U_{n,m}^{**}:=H(n,4;1^0,1^{n-2m}2^{m-2},1^0), &  U_{n,m}'':=H(n,4;1^1,1^{n-2m+1}2^{m-3},1^0),\\[5pt]
    U_{n,m}^*:=H(n,3;1^{n-2m+1}2^{m-3},1^1,1^1),& U_{n,m}^{3*}:=H(n,3;1^{n-2m}2^{m-2},1^1,1^0).
\end{array}
\end{equation*}
Obviously, $U_{n,m},\,U_{n,m}',\,U_{n,m}'',\, U_{n,m}^*,\,U_{n,m}^{**},\,U_{n,m}^{3*}\in \mathscr{U}^2_{n,m}.$

We denote by $H_{n,k}$ the graph obtained from $C_k$ by adding $n-k$ pendant vertices to a vertex of $C_k,$ and $U_n(k)$ the unicyclic graph obtained from $C_k$ by attaching one pendant vertex and $n-k-1$ pendant vertices to two adjacent vertices of $C_k$, respectively, where $3 \leq k \leq n-2.$ In \cite{12}, it was shown that
\begin{eqnarray}
  \xi^d(H_{n,k}) &=& \left\{
  \begin{array}{ll}
    -\frac{1}{8}k^4+\frac{n-1}{4}k^3+\frac{7-3n}{4}k^2+(n^2-\frac{7}{2}n+1)k+2n^2+n-n, & \hbox{if $k$ is even;} \\
    \\
    -\frac{1}{8}k^4+\frac{2n-1}{8}k^3+\frac{13-8n}{8}k^2+(n^2-\frac{9}{4}n+\frac{1}{8})k+n^2+\frac{1}{2}, & \hbox{if $k$ is odd.}
  \end{array}
\right.\label{eq:2.3} \\[3pt]
  \xi^d(U_n(k)) &=& \left\{
  \begin{array}{ll}
    -\frac{1}{8}k^4+\frac{n-1}{4}k^3+\frac{2-3n}{4}k^2+(n^2-2n-2)k+2n^2+n-2, & \hbox{if $k$ is even;} \\
    \\
    6n^2-11n-15, & \hbox{if $k=3$;} \\
    \\
   -\frac{1}{8}k^4+\frac{2n-1}{8}k^3+\frac{3-8n}{8}k^2+(n^2-\frac{3}{4}n-\frac{11}{8})k+n^2+\frac{1}{2}n-\frac{7}{4} , & \hbox{if $k\geq5$ is odd.}
  \end{array}
\right.\label{eq:2.2}
\end{eqnarray}
We know from \cite{12} that $H_{n,k}$ and $U_n(k)$ are, respectively, the extremal unicyclic graphs among $n$-vertex unicyclic graphs of girth $k$ having the minimal and second minimal eccentric distance sums.

The following formulas are easy to verify (see \cite{1}) that
\begin{equation}\label{eq:2.4}
W(C_k)=\left\{
  \begin{array}{ll}
    \frac{k^3}{8}, & \hbox{if $k$ is even;} \\[5pt]
   \frac{k(k^2-1)}{8}, & \hbox{if $k$ is odd,}
  \end{array}
\right.\ \ \
\text{and} \ \ \ \ D_{C_k}(v_i)=\left\{
  \begin{array}{ll}
    \frac{k^2}{4}, & \hbox{if $k$ is even;} \\[5pt]
   \frac{k^2-1}{4}, & \hbox{if $k$ is odd,}
  \end{array}
  \right. \text{for $v_i \in V_{C_k}.$}
\end{equation}

Further on we need the following lemmas.
\begin{lem}
[\cite{C-R-S1}]
Let $G \in \mathscr{U}_{2m,m}$ with $m \geq 3,$ and let $T$ be a pendant tree of $G$ with root $r$. If $u \in V_T$ is a pendant vertex furthest from the root $r$ with $d_G(u,r) \geq 2,$ then $u$ is adjacent to a vertex of degree two.
\end{lem}

\begin{lem}[\cite{D-I}]
A matching $M$ in $G$ is a maximum matching if and only if $G$ contains no $M$-augmenting path.
\end{lem}

\begin{lem}[\cite{2}]
Let $G \in \mathscr{U}_{n,m}$ with $n > 2m$ and $G \ncong C_n.$ Then there is a maximum matching $M$ and a pendant vertex $v$ of $G$ such that $v$ is not $M$-saturated.
\end{lem}

\begin{lem}
Given an $n$-vertex unicyclic graph $G$ with $\Delta(G)< n-1.$ Suppose that $u$ is a pendant vertex of $G$ whose neighbor is $v$, then
\[\label{eq:2.4}
\xi^d(G)\ge \xi^d(G-u) -3d_G(v)+9n-10+2\sum\limits_{x \in N(v)\backslash \{u\}} \varepsilon_G(x)+3\sum\limits_{x \in V_G\backslash N[v]}\varepsilon_G(x),
\]
with equality if and only if $\varepsilon_G(v)=2$ and $\varepsilon_G(x)=\varepsilon_{G-u}(x)$ for all $x \in V_{G-u}$.

Moreover, if $N(v)=\{u,w\},$ then
\[\label{eq:2.5}
\xi^d(G)\ge \xi^d(G-u-v)-7d_G(w)+25n-54+5\sum\limits_{x \in N[w] \backslash \{v\}}\varepsilon_G(x)+7\sum\limits_{x \in V_G\backslash (N[w]\bigcup \{u\})}\varepsilon_G(x),
\]
with equality if and only if $\varepsilon_G(w)=2$ and $\varepsilon_G(x)=\varepsilon_{G-u-v}(x)$ for all $x \in V_{G-u-v}$.
\end{lem}
\begin{proof}
By the definition of EDS, we have
\begin{eqnarray}
\xi^d(G)       &=& \sum\limits_{x \in V_{G-u}}\varepsilon_G(x)D_G(x) + \varepsilon_G(u)D_G(u)\notag\\
               &=& \sum\limits_{x \in V_{G-u}}\varepsilon_G(x)(D_{G-u}(x)+d_G(x,u))+\varepsilon_G(u)\sum\limits_{x \in V_G}d_G(x,u) \notag\\
               &=& \sum\limits_{x \in V_{G-u}}\varepsilon_G(x)D_{G-u}(x) +\sum\limits_{x \in V_{G-u}}\varepsilon_{G}(x)d_{G}(x,u) +\varepsilon_G(u)\sum\limits_{x \in V_G}d_G(x,u) \notag\\
               &=& \sum\limits_{x \in V_{G-u}}\varepsilon_{G-u}(x)D_{G-u}(x)+\sum\limits_{x \in V_{G-u}}(\varepsilon_G(x)-\varepsilon_{G-u}(x))D_{G-u}(x)\notag\\
               && +\varepsilon_G(u)\sum\limits_{x \in V_G}d_G(x,u)+\sum\limits_{x \in V_{G-u}}\varepsilon_{G}(x)d_{G}(x,u) \notag\\
               &\geq& \xi^d(G-u) + \Theta+\Theta',\label{2.6}
\end{eqnarray}
where
$$
\Theta = \varepsilon_G(u)\sum\limits_{x \in V_G}d_G(x,u),\ \
\Theta'= \sum\limits_{x \in V_{G-u}}\varepsilon_{G}(x)d_{G}(x,u).
$$
The last inequality follows from $\varepsilon_G(x) \geq \varepsilon_{G-u}(x)$. The equality in (\ref{2.6}) holds if and only if $\varepsilon_G(x)=\varepsilon_{G-u}(x)$ for all $x \in V_{G-u}.$ In what follows we are to determine sharp lower bounds on $\Theta$ and $\Theta'$, respectively.

On the one hand,
\begin{eqnarray}
\Theta & =& \varepsilon_G(u)\sum_{x \in V_G}d_G(x,u)\notag  \\
      & =& \varepsilon_G(u)\left(\sum\limits_{x \in N[v] \backslash \{u\}}d_G(x,u)+ \sum\limits_{x \in V_G \backslash N[v]}d_G(x,u)\right) \notag\\
      & =& \varepsilon_G(u)\left(1+2(d_G(v)-1)+\sum\limits_{x \in V_G \backslash N[v]}d_G(x,u)\right)\notag  \\
      & \geq & 3 (1+2(d_G(v)-1)+3(n-d_G(v)-1)) \label{2.7}\\
      & = & 9n-3d_G(v)-12.\notag
\end{eqnarray}
Note that $\Delta(G)< n-1$, hence $V_G \backslash N[v]\not=\emptyset.$ The equality in (\ref{2.7}) holds if and only if $\varepsilon_G(u)=3$ and $d_G(x,u)=3$ for all $x \in V_{G-u} \backslash N[v]$, which is equivalent to $\varepsilon_G(v)=2$.

On the other hand,
\begin{eqnarray}
\Theta'& =& \sum\limits_{x \in V_{G-u}}\varepsilon_{G}(x)d_{G}(x,u)  \notag\\
      & =& \varepsilon_G(v)d_G(v,u)+ \sum\limits_{x \in N(v)\backslash \{u\}}\varepsilon_G(x)d_G(x,u) +\sum\limits_{x \in V_{G-u} \backslash N[v]}\varepsilon_G(x)d_G(x,u) \notag \\
      & = & 2+\sum\limits_{x \in N[v]\backslash \{u,v\}}2\varepsilon_G(x)+ \sum\limits_{x \in V_{G-u} \backslash N[v]}\varepsilon_G(x)d_G(x,u) \notag \\
      & \geq & 2+\sum\limits_{x \in N[v]\backslash \{u,v\}}2\varepsilon_G(x)+ \sum\limits_{x \in V_{G-u} \backslash N[v]}3\varepsilon_G(x)\label{2.8}
\end{eqnarray}
Note that $\Delta(G)< n-1$, hence $V_{G-u} \backslash N[v]\not=\emptyset.$ The equality in (\ref{2.8}) holds if and only if $d_G(x,u)=3$ for all $x \in V_{G-u} \backslash N[v]$, which is equivalent to $\varepsilon_G(v)=2$.

Combining with (\ref{2.6})-(\ref{2.8}) we obtain inequality in (\ref{eq:2.4}) holds, with equality if and only if $\varepsilon_G(v)=2$ and $\varepsilon_G(x)=\varepsilon_{G-u}(x)$ for all $x \in V_{G-u}$. This completes the proof of the first part of Lemma 2.4.

Note that $u$ is a pendant and $N(v)=\{u,w\},$ hence we have
\begin{eqnarray}
\xi^d(G)       &=& \sum\limits_{x \in V_{G-u-v}}\varepsilon_G(x)D_G(x) + \varepsilon_G(u)D_G(u)+\varepsilon_G(v)D_G(v)\notag\\
               &=& \sum\limits_{x \in V_{G-u-v}}\varepsilon_G(x)[D_{G-u-v}(x)+(d(x,w)+2)+(d(x,w)+1)]+\varepsilon_G(u)D_G(u)+\varepsilon_G(v)D_G(v) \notag\\
               &=& \sum\limits_{x \in V_{G-u-v}}\varepsilon_G(x)D_{G-u-v}(x) + \sum\limits_{x \in V_{G-u-v}}\varepsilon_G(x)(2d_G(x,w)+3) +\varepsilon_G(u)D_G(u)+\varepsilon_G(v)D_G(v) \notag\\
               &=& \sum\limits_{x \in V_{G-u-v}}\varepsilon_{G-u-v}(x)D_{G-u-v}(x)+\sum\limits_{x \in V_{G-u-v}}(\varepsilon_G(x)-\varepsilon_{G-u-v}(x))D_{G-u-v}(x)  \notag\\
               &&+\varepsilon_G(u)D_G(u)+\varepsilon_G(v)D_G(v)+\sum_{x \in V_{G-u-v}}\varepsilon_G(x)(2d_G(x,w)+3) \notag\\
               &\geq& \xi^d(G-u-v) + \Theta_1 +\Theta_1', \label{2.9}
\end{eqnarray}
where
\begin{eqnarray*}
\Theta_1 = \varepsilon_G(u)D_G(u)+\varepsilon_G(v)D_G(v), \ \ \ \ \Theta_1' =\sum_{x \in V_{G-u-v}}\varepsilon_G(x)(2d_G(x,w)+3).
\end{eqnarray*}
The last inequality follows from $\varepsilon_G(x) \geq \varepsilon_{G-u-v}(x)$. The equality in (\ref{2.9}) holds if and only if $\varepsilon_G(x)=\varepsilon_{G-u-v}(x)$ for all $x \in V_{G-u-v}.$ In what follows we are to determine sharp lower bounds on $\Theta_1$ and $\Theta_1'$, respectively.

On the one hand,
\begin{eqnarray}
\Theta_1 & =& \varepsilon_G(u)D_G(u)+\varepsilon_G(v)D_G(v) \notag \\
      & = & (\varepsilon_G(w)+2)(D_G(w)+2n-6)+(\varepsilon_G(w)+1)(D_G(w)+n-4)  \notag \\
      & = &(\varepsilon_G(w)+2)\left(\sum\limits_{x \in N(w)}d_G(x,w)+\sum\limits_{x \in V_G\backslash N[w]}d_G(x,w)+2n-6\right) \notag \\
      & &  +(\varepsilon_G(w)+1)\left(\sum\limits_{x \in N(w)}d_G(x,w)+\sum\limits_{x \in V_G\backslash N[w]}d_G(x,w)+n-4\right)  \notag \\
      & \geq & 4(d_G(w)+2(n-d_G(w)-1)+2n-6)+3(d_G(w)+2(n-d_G(w)-1)+n-4)\label{2.10}\\
      & =& 25n-7d_G(w)-50. \notag
\end{eqnarray}
Note that $\Delta(G)< n-1$, hence $V_G \backslash N[v]\not=\emptyset.$ The equality in (\ref{2.10}) holds if and only if $\varepsilon_G(w)=2$ and $d(x,w)=2$ for all $x\in V_G\setminus N[w]$, which is equivalent to that $\varepsilon_G(w)=2$.

On the other hand,
\begin{eqnarray}
\Theta_1' 
      & = & 2\sum\limits_{x \in V_{G-u-v}}\varepsilon_G(x)d_G(x,w)+3\sum\limits_{x \in V_{G-u-v}}\varepsilon_G(x) \notag \\
      & = & 2\left[\sum\limits_{x \in N[w] \backslash \{w,v\}}\varepsilon_G(x)+ \sum\limits_{x \in V_{G-u-v}\backslash N[w]}\varepsilon_G(x)d_G(x,w)\right]+3\sum\limits_{x \in V_{G-u-v}}\varepsilon_G(x)\notag \\
      & \geq & 2\left[\sum\limits_{x \in N(w) \backslash \{v\}}\varepsilon_G(x)+ 2\sum\limits_{x \in V_{G-u-v}\backslash N[w]}\varepsilon_G(x)\right]+3\sum\limits_{x \in V_{G-u-v}}\varepsilon_G(x)\label{2.11}\\
      & = & 2\left[\sum\limits_{x \in V_{G-u-v} \backslash \{w\}}\varepsilon_G(x)+ \sum\limits_{x \in V_{G-u-v}\backslash N[w]}\varepsilon_G(x)\right]+3\sum\limits_{x \in V_{G-u-v}}\varepsilon_G(x) \notag \\
      & = & 5\sum\limits_{x \in N[w] \backslash \{v\}}\varepsilon_G(x)+7\sum\limits_{x \in V_G\backslash \{N[w]\bigcup \{u\}\}}\varepsilon_G(x)-4. \notag
\end{eqnarray}
The equality in (\ref{2.11}) holds if and only if $d(x,w)=2$ for all $x \in V_{G-u-v}\setminus N[w]\not=\emptyset$, which is equivalent to that $\varepsilon_G(w)=2$ and $V_{G-u-v}\setminus N[w]\not=\emptyset.$

Combining with (\ref{2.9})-(\ref{2.11}) yields inequality (\ref{eq:2.5}), with equality if and only if $\varepsilon_G(w)=2$ and $\varepsilon_G(x)=\varepsilon_{G-u-v}(x)$ for all $x \in V_{G-u-v}$ and $V_{G-u-v}\setminus N[w]\not=\emptyset,$ which is equivalent to
$\varepsilon_G(w)=2$ and $\varepsilon_G(x)=\varepsilon_{G-u-v}(x)$ for all $x \in V_{G-u-v}.$

This completes the proof of the second part of Lemma 2.4.
\end{proof}

\begin{lem}
Suppose that $m+1 \leq k \leq 2m-2$ with $m \geq 5$, then $\xi^d(U_{2m}(k)) > 43m^2 -72m + 6.$
\end{lem}
\begin{proof}
We show our result according to the parity of $k$ with $k \geq 6.$ Here we only show our result holds for even $k$. Similarly, we can also show that our result holds for odd $k$, which is omitted here.

For even $k$, in view of (\ref{eq:2.2}) we have
$$
\xi^d(U_{2m}(k))=-\frac{1}{8}k^4+\frac{2m-1}{4}k^3+\frac{2-6m}{4}k^2+(4m^2-4m-2)k+8m^2+2m-2.
$$

Let
$$
f(x)=-\frac{1}{8}x^4+\frac{2m-1}{4}x^3+\frac{2-6m}{4}x^2+(4m^2-4m-2)x+8m^2+2m-2,
$$
where $6\le m+1 \leq x \leq 2m-2.$
It is easy to check that
$$
f'(x)=-\frac{1}{2}x^3+\frac{3}{4}(2m-1)x^2+(1-3m)x+4m^2-4m-2,\, \, f''(k)=-\frac{1}{2}(3x^2-3(2m-1)x+6m-2).
$$
Let $\varphi(x)=3x^2-3(2m-1)x+6m-2$, then we have $\varphi'(x)=3(2x-2m+1)>0 \,\,\textrm{for}\,\,  m+1 \leq x \leq 2m-2.$
$$
\varphi(m+1)=-6m^2+6m+4 < 0,\,\,\, \varphi(2m-2)=4 > 0.
$$
Therefore, there exists a $c_0 \in (m+1,2m-2)$, $f''(x)$ is positive in $[m+1,c_0)$, negative in $(c_0,2m-2]$. Hence $f'(x)$ is increasing in $[m+1,c_0)$ and decreasing in $(c_0,2m-2]$. Thus, for $x \in [m+1,2m-2],$ $f'(x)$ takes its minimal value at $k=m+1$ or $k=2m-2$. On the other hand, for $m\ge 5$,
$$
f'(m+1)= (m+1)^2(m-\frac{5}{4})+(m-3)^2-10>0,\ \ \
f'(2m-2)= 2(m-1)^3+(m-1)^2-2>0.
$$
So we get that $f'(x)>0$ for $x \in [m+1,2m-2],$ which implies that $f(x)$ is an increasing function in $[m+1,2m-2].$ Note that
 \begin{align*}
 f(m+1)-(43m^2-72m+6)
                       =& \frac{3}{8}(m+1)^4+\frac{7}{4}(m+1)^3-45(m+1)^2+150m+33.
\end{align*}
It is easy to check that $\Phi(m):=\frac{3}{8}(m+1)^4+\frac{7}{4}(m+1)^3-45(m+1)^2+150m+33$ is increasing in $[5,+\infty)$, so $\Phi(m)\geq \Phi(5)=27 >0.$

This completes the proof.
\end{proof}

For integer $m \geqslant 3,$ let $\mathscr{U}'_m$ be the set of graphs in $\mathscr{U}_{2m,m}$ containing a pendant vertex whose neighbor is of degree two. Let $\mathscr{U}''_m=\mathscr{U}_{2m,m} \backslash \mathscr{U}'_m.$

\begin{lem}
Let $G \in \mathscr{U}''_m$ with $m \geq 5$. Then $\xi^d(G)> 43m^2-72m+6.$
\end{lem}
\begin{proof}
If $G \cong C_{2m},$ $\xi^d(C_{2m})=2m^4 > 43m^2-72m+6,$ as desired.

Suppose that $G \ncong C_{2m}.$ By Lemma 2.1, $G$ is a graph of maximum degree three obtained by attaching some pendant vertices to a cycle $C_k,$ where $m \leq k \leq 2m-1.$

If $k=m,$ then every vertex on the cycle has degree three. So it is easy to check that
\begin{equation*}
\xi^d(G)=\left\{
  \begin{array}{ll}
    \sum_{i=1}^m[(\frac{m}{2}+1)(2D_{C_{m}}(v)+m)+(\frac{m}{2}+2)(2D_{C_{m}}(v)+3m-2)], & \hbox{if $m$ is even;} \\[5pt]
   \sum_{i=1}^m[\frac{m+1}{2}(2D_{C_{m}}(v)+m)+(\frac{m+1}{2}+1)(2D_{C_{m}}(v)+3m-2)], & \hbox{if $m$ is odd.}
  \end{array}
\right.
\end{equation*}
Together with (\ref{eq:2.4}), if $k=m \geq 5$ is even, then
$$
\xi^d(G)=\frac{1}{2}m^4+\frac{7}{2}m^3+6m^2-4m > 43m^2-72m+6;
$$
if $k=m \geq 5 $ is odd, then we have
$$
\xi^d(G)=\frac{1}{2}m^4+3m^3+\frac{7}{2}m^2-4m  > 43m^2-72m+6.
$$

If $m+1 \leq k \leq 2m-2,$ then for a given $k\ge 6$, we know from \cite{12} that $H_{n,k}, U_{2m}(k)$ are, respectively, the graphs with the minimal and second minimal EDS among $n$-vertex unicyclic graphs of girth $k.$ It is easy to see that $G\not\cong H_{n,k}$ (since $H_{n,k}$ contains no perfect matching). Together with Lemma 2.5, for some $U_{2m}(k),$ we have $\xi^d(G) \ge \xi^d(U_{2m}(k)) > 43m^2-72m+6.$

If $k=2m-1,$ then $G$ is the graph obtained from $C_k$ by attaching just one pendant vertex to a vertex of $C_k$. By direct computation, we have $\xi^d(G)=2m^4-3m^3+7m^2-4m+1 > 43m^2-72m+6$ for $m \geq 5.$

This completes the proof.
\end{proof}

It is routine to check that if $G \in \mathscr{U}_{n,m},$ then
\[\label{eq:2.12}
\Delta(G)\leq n-m+1
 \]
with equality if and only if $G \cong U_{n,m}$. Furthermore, if $G \in \mathscr{U}_{n,m}\backslash\{U_{n,m}\},$ then
\[\label{eq:2.13}
\Delta(G)\leq n-m
\]
with equality if and only if $G$ is isomorphic to one of graphs in the set $\{U_{n,m}', U_{n,m}'', U_{n,m}^*, U_{n,m}^{**}, U_{n,m}^{3*}\}.$

\section{\normalsize  Lower bounds and the extremal graphs }\setcounter{equation}{0}
In this section, we are to determine the first and second minimal EDS of graphs in $\mathscr{U}_{n,m}$, the corresponding extremal graphs are identified.
\begin{table}[h!]
  \centering
  \caption{The minimal EDS and the corresponding extremal graphs in $\mathscr{U}_{n,m}$ with small order}\label{dd}\vspace{2mm}

    \begin{tabular}{llllll}
                         \hline
                         $n \backslash m$ & 2 & 3 & 4 & 5 & 6\\ \hline
                         4  &  $29,\,  U_{4,2}$ &     &     &     & \\[1pt]
                         5  & $54,\,  U_{5,2}$  &     &     &     & \\[1pt]
                         6  & $91,\,  U_{6,2}$  & $133,\,  U_{6,3}'$ &     &     & \\[1pt]
                         7  & $134,\,  U_{7,2}$ & $206,\,  U_{7,3}'$ &     &     & \\[1pt]
                         8  & $185,\,  U_{8,2}$ & $291,\,  U_{8,3}'$ & $373,\,  C_5(1^1,1^1,1^1)$ &     & \\[1pt]
                         9  & $244,\,  U_{9,2}$ & $388,\,  U_{9,3}'$ & $484,\,  U_{9,4}$ &     & \\[1pt]
                         10 & $311,\,  U_{10,2}$ & $496,\,  U_{10,3}$ & $603,\,  U_{10,4}$ & $672,\,  U_{10,5}$ & \\[1pt]
                         11 & $386,\,  U_{11,2}$ & $613,\,  U_{11,3}$ & $734,\,  U_{11,4}$ & $812,\,  U_{11,5}$ & \\[1pt]
                         12 & $469,\,  U_{12,2}$ & $742,\,  U_{12,3}$ & $877,\,  U_{12,4}$ & $964,\,  U_{12,5}$ & $1053,\,  U_{12,6}$\\
                         \hline
\end{tabular}
\end{table}

For small $n$, with the help of Nauty \cite{19} we may list all the possible unicyclic graphs in $\mathscr{U}_{n,m}$ for $4\leqslant n\leqslant 16, 2\leqslant m\leqslant 6$.
Table 1 presents the smallest values of EDS and the corresponding extremal graphs in $\mathscr{U}_{n,m}$ with $4\leqslant n\leqslant 12, 2\leqslant m\leqslant 6$. In Table 2, we list the second smallest values of EDS and the corresponding extremal graphs in $\mathscr{U}_{n,m},$ where $4\leqslant n\leqslant 12, 2\leqslant m\leqslant 6$ and $m=4, n=13,14,15,16$.

We first determine the first and second smallest EDS of unicyclic graphs with perfect matchings. It is routine to check that $\mathscr{U}_{4,2}=\{U_{4,2}, C_4\}$ and $\mathscr{U}_{6,3}=\{U_{6,3}, U_{6,3}', U_{6,3}'',U_{6,3}^*, U_{6,3}^{**}, U_{6,3}^{3*},C_6\}.$ By direct calculation, we know that, for $m=2,3$,
$$
\xi^d(U_{4,2})<\xi^d(C_4)
$$
and
$$
\xi^d(U_{6,3}')<\xi^d(U_{6,3}^*)<\xi^d(U_{6,3})=\xi^d(U_{6,3}'')<\xi^d(C_6)<\xi^d(U_{6,3}^{**})<\xi^d(U_{6,3}^{3*}).
$$
Hence, in what follows we consider $m\ge 4$ for $\mathscr{U}_{2m,m}.$
\begin{table}[h!]
  \centering
  \caption{The second minimal EDS and the corresponding extremal graphs in $\mathscr{U}_{n,m}$ with small order}\label{dd}\vspace{2mm}

    \begin{tabular}{llllll}
                         \hline
                         $n \backslash m$ & 2 & 3 & 4 & 5 & 6\\ \hline
                         4  &  $32,\,  C_4$ &     &     &     & \\[1pt]
                         5  & $60,\,  C_5$  &     &     &     & \\[1pt]
                         6  & $134,\,  U_{6,2}^{**}$  & $141,\,  U_{6,3}^*$ &     &     & \\[1pt]
                         7  & $201,\,  U_{7,2}^{**}$ & $214,\,  U_{7,3}^*$ &     &     & \\[1pt]
                         8  & $280,\,  U_{8,2}^{**}$ & $298,\,  U_{8,3}$ & $377,\,  U_{8,4}$ &     & \\[1pt]
                         9  & $371,\,  U_{9,2}^{**}$ & $391,\,  U_{9,3}$ & $492,\,  C_5(1^1,1^2,1^1)$ &     & \\[1pt]
                         10 & $474,\,  U_{10,2}^{**}$ & $497,\,  U_{10,3}'$ & $623,\,  C_5(1^1,1^3,1^1)$ & $711,\,  U_{10,5}'$ & \\[1pt]
                         11 & $589,\,  U_{11,2}^{**}$ & $618,\,  U_{11,3}'$ & $766,\,  C_5(1^1,1^4,1^1)$ & $860,\,  U_{11,5}'$ & \\[1pt]
                         12 & $716,\,  U_{12,2}^{**}$ & $751,\,  U_{12,3}'$ & $921,\,  C_5(1^1,1^5,1^1)$ & $1021,\,  U_{12,5}'$ & $1112,\,  U_{12,6}'$\\[1pt]
                         13  & -----                       &----- & $1088, \,C_5(1^1,1^6,1^1)$ & ----- & -----\\[1pt]
                         14  & -----                       &----- & $1267, \,C_5(1^1,1^7,1^1)$ & ----- & -----\\[1pt]
                         15  & -----                       &----- & $1458, \,C_5(1^1,1^8,1^1)$ & ----- & -----\\[1pt]
                         16  & -----                       &----- & $1660,\,U_{16,4}'$ & ----- & -----\\
                         \hline
\end{tabular}
\end{table}
\begin{thm}
Let $G \in \mathscr{U}_{2m,m}$ with $m \geq 4.$
\begin{wst}
\item[{\rm (i)}]If $m=4$, then $\xi^d(G)\geq 373$ with equality if and only if $G \cong C_5(1^1,1^1,1^1)$.
\item[{\rm (ii)}]If $m \geq 5$, then $\xi^d(G) \geq  43m^2 - 92m + 57$ with equality if and only if $G \cong U_{2m,m}$.
\end{wst}
\end{thm}
\begin{proof}
(i)\ \ The cases for $m=4$, our result follows directly from Table 1.

(ii)\ \ Suppose that $m \geq 5$ and let $g_1(m)=43m^2-92m+57$ in what follows. We prove the result by induction on $m.$

If $m=5,$ then the result follows from Table 1. Suppose that $m \geq 6$ and the result holds for graphs in $\mathscr{U}_{2k-2,k-1}$ with $k\le m-1$. Now we consider $k=m$ and let $G \in \mathscr{U}_{2m,m}$.

If $G \in \mathscr{U}''_m,$ then by Lemma~2.6, $\xi^d(G) > 43m^2-72m+6 > g_1(m).$ If $G \in \mathscr{U}_m',$ then $G$ contains a pendant vertex, say $u$, such that its unique neighbor, say $v$, is of degree 2. Hence, $G-u-v \in \mathscr{U}_{2m-2,m-1}.$ For convenience, let $N_G(v)=\{u,w\}$. By first part of Lemma 2.4, it is easily seen that
\begin{eqnarray}
 \xi^d(G)&\geq& \xi^d(G-u-v)-7d_G(w)+25n-54+\sum_{x \in N[w] \backslash \{v\}}5\varepsilon_G(x)+\sum_{x \in V_G\backslash \{N[w]\bigcup \{u\}\}}7\varepsilon_G(x)\label{eq:3.1}\\
         &=&  \xi^d(G-u-v)-20d_G(w)+106m-115 \\
         & \geq& g_1(m-1)-20(m+1)+106m-115 \label{eq:3.2}\text{\ \ \ \ \ \ (by the induction hypothesis and (\ref{eq:2.12}))}\\
         &=& g_1(m). \notag
\end{eqnarray}
Based on Lemma 2.4, the equality in (\ref{eq:3.1}) holds if and only if $\varepsilon_G(w)=2$ and $\varepsilon_{G-u-v}(x)=\varepsilon_G(x)$ for all $x \in V_{G-u-v},$ which (based on $m\ge 6$) implies that $G-u-v$ is in
$$
   \left\{H(2m-2,5;1^1,1^12^{m-4},1^1),U_{2m-2,m-1},U_{2m-2,m-1}',U_{2m-2,m-1}'',U_{2m-2,m-1}^*,U_{2m-2,m-1}^{**},U_{2m-2,m-1}^{3*}\right\}.
$$
Hence, we obtain that $\varepsilon_G(x)=3$ for all $x \in N(w)$ and $\varepsilon_G(x)=4$ for $x \in V_G\backslash N[w]$. Hence equality in (3.2) holds.
The equality holds in (\ref{eq:3.2}) if and only if $G-u-v \cong U_{2m-2,m-1}$, $d_G(w)=m+1.$ That is to say, $\xi^d(G)=g_1(m)$ holds if and only if $G\cong U_{2m,m}.$

This completes the proof.
\end{proof}

From Table 2, it is obvious that $U_{8,4}$ is the unique graph with second minimal EDS in $\mathscr{U}_{8,4}.$ Then we consider $m \geq 5$ for $\mathscr{U}_{2m,m}.$
\begin{thm}
Let $G \in \mathscr{U}_{2m,m}\backslash\{U_{2m,m}\}$ with $m \geq 5.$ Then $\xi^d(G) \geq  43m^2 - 72m -4$ with equality if and only if $G \cong U_{2m,m}'$.
\end{thm}

\begin{proof}Suppose that $m \geq 5$ and let $g_2(m)=43m^2-72m-4$. We are proceed by induction on $m.$

If $m=5,$ then the result follows from Table 2. Suppose that $m \geq 6$ and the result holds for graphs in $\mathscr{U}_{2k-2,k-1}$ with $k\le m-1$. Now we consider $k=m$. Let $G \in \mathscr{U}_{2m,m}\backslash\{U_{2m,m}\}$. If $G \in \mathscr{U}''_m,$ then by Lemma~2.6, $\xi^d(G) > 43m^2-72m+6 > g_2(m).$ If $G \in \mathscr{U}_m',$ then $G$ contains a pendant vertex, say $u$, such that its unique neighbor, say $v$, is of degree 2. Hence, $G-u-v \in \mathscr{U}_{2m-2,m-1}\backslash\{U_{2m-2,m-1}\}.$ For convenience, let $N_G(v)=\{u,w\}$. By the second part of Lemma 2.4, we have
\begin{eqnarray}
 \xi^d(G)&\geq& \xi^d(G-u-v)-7d_G(w)+25n-54+5\sum\limits_{x \in N[w] \backslash \{v\}}\varepsilon_G(x)+7\sum\limits_{x \in V_G\backslash \{N[w]\bigcup \{u\}\}}\varepsilon_G(x)\label{3.4}\\
         &=&  \xi^d(G-u-v)-20d_G(w)+106m-115 \\
         & \geq& g_2(m-1)-20m+106m-115 \label{3.5}\text{\ \ \ \ \ \ (by the induction hypothesis and (\ref{eq:2.13}))}\\
         &=& g_2(m). \notag
\end{eqnarray}
Based on Lemma 2.4, the equality in (\ref{3.4}) holds if and only if $\varepsilon_G(w)=2$ and $\varepsilon_{G-u-v}(x)=\varepsilon_G(x)$ for all $x \in V_{G-u-v},$ which (based on $m\ge 6$) implies that
$$
  G-u-v \in \{H(2m-2,5;1^1,1^12^{m-5},1^1),U_{2m-2,m-1}',U_{2m-2,m-1}'',U_{2m-2,m-1}^*,U_{2m-2,m-1}^{**},U_{2m-2,m-1}^{3*}\}.
$$
Hence, we obtain that $\varepsilon_G(x)=3$ for all $x \in N(w)$, $\varepsilon_G(x)=4$ for $x \in V_G\backslash N[w]$, which implies that equality in (3.5) holds. The equality in (\ref{3.5}) holds if and only if $G-u-v \cong U_{2m-2,m-1}'$, $d_G(w)=m.$ That is to say, $\xi^d(G)=g_2(m)$ holds if and only if $G\cong U_{2m,m}'.$

This completes the proof.
\end{proof}

In the rest of this section, we are to determine the graphs in $\mathscr{U}_{n,m}$ with the first and second minimal EDS with $m\ge 2.$ Note that for $\mathscr{U}_{n,3}$ with $n\le 9$, we can easily determine the graph with minimal EDS from Table 1.
\begin{thm}
Let $G$ be an $n$-vertex unicyclic graph with matching number $m$, where $2 \leq m \leq \lfloor\frac{n}{2}\rfloor.$ 
\begin{wst}
\item[{\rm (i)}]Among $\mathscr{U}_{n,2},$ one has $\xi^d(G)\geq 4n^2-9n+1$ with equality if and only if $G \cong U_{n,2}.$
\item[{\rm (ii)}]Among $\mathscr{U}_{n,3}$ with $n\le 9$, one has $\xi^d(G)\geq 6n^2-5n-53$ with equality if and only if $G \cong U_{n,3}';$ while among
$\mathscr{U}_{n,3}$ with $n\ge 10$, one has $\xi^d(G)\geq 6n^2-9n-14$ with equality if and only if $G \cong U_{n,3}.$
\item[{\rm (iii)}]Among $\mathscr{U}_{n,m}$ with $n\ge 9, m \geq 4,$ one has $\xi^d(G) \geq 6n^2 + m^2 +9mn - 30m - 31n +57$ with equality if and only if $G \cong U_{n,m}$.
\end{wst}
\end{thm}
\begin{proof}
(i)\ \ For the case $n=5,6,$ our result follows immediately from Table 1. Hence, in what follows we consider $n \geq 7.$ It is easy to see that $\mathscr{U}_{n,2}=\{C_3(1^{k_1},1^{k_2},1^0)\}\cup \{H_{n,4}\}\cup \{H(n,3;1^0S_{n-4},1^0,1^0)\}.$

Note that $U_{n,2}\in \{C_3(1^{k_1},1^{k_2},1^0):k_1+k_2=n-3\},$ and from \cite{11} we know that $U_{n,2}$ is the graph with the minimal EDS among $n$-vertex unicyclic graphs, and of course it is the graph with minimal degree distance in $\mathscr{U}_{n,2}$, as desired.

(ii)\ \  In fact, among $\mathscr{U}_{n,3}$ with $n\le 9$, our result follows directly from Table 1. So we consider $m=3, n\ge 10$ in what follows. In this case, it is easy to see that $G\not\cong C_n.$ We show our result by induction on $n.$ If $n=10$, then our result holds immediately from Table 1. So we consider that $n \geq 11$ and assume that our result holds for graphs in $\mathscr{U}_{k,3}$ with $k \leq n-1.$ In view of Lemma 2.4, we have
inequality in (\ref{eq:2.4}) holds,
with equality if and only if $\varepsilon_G(v)=2$ and $\varepsilon_{G-u}(x)=\varepsilon_G(x)$ for all $x \in V_{G-u},$ where $u$ is a pendant vertex and $v$ is its unique neighbor. Hence, by the definition of $\mathscr{U}^2_{n,3},$ we have
\begin{eqnarray*}
&&G\in \mathscr{U}^2_{n,3}= \left\{H(n,5;1^{k_1},1^{k_2+1}2^0,1^0)\right\}\bigcup \left\{H(n,4;1^{k_1+1},1^{k_2+1}2^0,1^{k_3})\right\}\bigcup \left\{H(n,4;1^{k_1},1^{k_2}S_{t_1},1^{k_3})\right\}\\
                 & & \ \ \ \ \ \ \ \ \ \ \ \ \  \,\,\,\,\,\,\,\, \bigcup\left\{H(n,3;1^{k_1+1}2^0,1^{k_2+1},1^{k_3+1})\right\}\bigcup \left\{H(n,3;1^{k_1}S_{t_1},1^{k_2},1^0)\right\},\ \ \  k_1,k_2,k_3\geq 0, t_1 \geq 1.
\end{eqnarray*}
In view of (\ref{eq:2.4}), let $\phi(G)=2\sum_{x \in N(v)\backslash \{u\}} \varepsilon_G(x)+3\sum_{x \in V_G\backslash N[v]}\varepsilon_G(x).$ By the structure of graphs in $\mathscr{U}^2_{n,3}$ and by direct calculation, we have
\begin{align*}
\phi(H(n,5;1^{k_1},1^{k_2+1}2^0,1^0))\geq &\ \phi(H(n,5;1^0,1^{n-5},1^0))=6n-10,\\
\phi(H(n,4;1^{k_1+1},1^{k_2+1}2^0,1^{k_3})) \geq& \ \phi(H(n,4;1^1,1^{n-5},1^0))=6n-8 ,\\
\phi(H(n,4;1^{k_1},1^{k_2}S_{t_1},1^{k_3}))\geq& \ \phi(H(n,4;1^0,1^{n-6}2^1,1^0))=6n,\\
\phi(H(n,3;1^{k_1+1}2^0,1^{k_2+1},1^{k_3+1}))\geq& \ \phi(H(n,3;1^{n-5}2^0,1^1,1^1))=6n-10,\\
\phi(H(n,3;1_{k_1}S_{t_1},1^{k_2},1^0))\geq & \ \phi(H(n,3;1^{n-5}2^1,1^0,1^0))=6n-11.
\end{align*}
Hence, in view of (\ref{eq:2.4}) and by induction we have
\begin{eqnarray}
\text{RHS of (\ref{eq:2.4})} &\geq& g(n-1,3)-3(n-2)+9n-10+6n-11 \ \ \ \ \text{(by induction and (\ref{eq:2.12}))} \label{3.6}\\
         &=& 6n^2-9n-14. \notag
\end{eqnarray}
The equality holds in (\ref{3.6}) if and only if $G-u \cong U_{n-1,3}$, $d_G(v)=n-2$, $ \sum_{x \in N(v)\backslash \{u\}} 2\varepsilon_G(x)+\sum_{x \in V_G\backslash N[v]}3\varepsilon_G(x) = 6n -11.$ That is to say, $\xi^d(G)=6n^2-9n-14$ holds if and only if $G \cong U_{n,3},$ as desired.

(iii)\ \ For convenience, let $f_1(n,m)=6n^2 + m^2 +9mn - 30m - 31n +57.$ If $G\cong C_n$, then $n=2m+1.$ In view of (2.3), we obtain that $\xi^d(C_{2m+1})=2m^4+3m^3+m^2.$ On the other hand,
$f_1(2m+1,m)=43m^2-59m+32.$ It is easy to check that $2m^4+3m^3+m^2>43m^2-59m+32.$ So in what follows we consider that $G\not\cong C_n.$

\medskip
\noindent$\bullet$ $m=4$. For $n=9$, our result follows directly from Table 1. Suppose that $n \geq 10$ and the result holds for graphs in $\mathscr{U}_{n-1,4}$. Hence, for graph $G\in \mathscr{U}_{n,4}$, in view of Lemma 2.4, we also have inequality holds in (\ref{eq:2.4}), with equality if and only if $\varepsilon_G(v)=2$ and $\varepsilon_{G-u}(x)=\varepsilon_G(x)$ for all $x \in V_{G-u},$ where $u$ is a pendant vertex and $v$ is its unique neighbor. Hence, by the definition of $\mathscr{U}^2_{n,4},$ we have
\begin{eqnarray*}
&&G\in \mathscr{U}^2_{n,4}= \left\{H(n,5;1^{k_1+1},1^{k_2+1}2^0,1^{k_3+1})\right\}\bigcup \left\{H(n,5;1^{k_1},1^{k_2}S_{t_1},1^{k_3})\right\}\bigcup \left\{H(n,4;1^{k_1},1^{k_2}S_{t_1,t_2},1^{k_3})\right\}\\
                 & & \ \ \ \ \ \ \ \ \ \ \ \ \  \,\,\,\,\,\,\,\, \bigcup\left\{H(n,3;1^{k_1+1}S_{t_1},1^{k_2+1},1^{k_3+1})\right\}\bigcup \left\{H(n,3;1^{k_1+1}S_{t_1,t_2},1^{k_2},1^{k_3})\right\},\   k_1,k_2,k_3\geq 0, t_1,t_2 \geq 1.
\end{eqnarray*}
In view of (\ref{eq:2.4}), let $\phi(G)=2\sum_{x \in N(v)\backslash \{u\}} \varepsilon_G(x)+3\sum_{x \in V_G\backslash N[v]}\varepsilon_G(x).$ By the structure of graphs in $\mathscr{U}^2_{n,4}$ and by direct calculation, we have
\begin{align*}
\phi(H(n,5;1^{k_1+1},1^{k_2+1}2^0,1^{k_3+1}))\geq&\ \phi(H(n,5;1^1,1^{n-7}2^0,1^1))=6n+6,\\
\phi(H(n,5;1^{k_1},1^{k_2}S_{t_1},1^{k_3})) \geq&\ \phi(H(n,5;1^0,1^{n-7}2^1,1^0))=6n+6,\\
\phi(H(n,4;1^{k_1+1},1^{k_2+1}S_{t_1},1^{k_3}))\geq&\ \phi(H(n,4;1^1,1^{n-7}2^1,1^0))=6n+6,\\
\phi(H(n,4;1^{k_1},1^{k_2}S_{t_1,t_2},1^{k_3}))\geq&\ \phi(H(n,4;1^0,1^{n-8}2^2,1^0))=6n+6,\\
\phi(H(n,3;1^{k_1+1}S_{t_1},1^{k_2+1},1^{k_3+1}))\geq&\ \phi(H(n,3;1^{n-7}2^1,1^1,1^1))=6n+6,\\
\phi(H(n,3;1^{k_1+1}S_{t_1,t_2},1^{k_2},1^{k_3}))\geq&\ \phi(H(n,3;1^{n-7}2^2,1^0,1^0))=6n.
\end{align*}
Hence, in view of (\ref{eq:2.4}) and by induction we have
\begin{eqnarray}
\text{RHS of (\ref{eq:2.4})} &\geq& f_1(n-1,4)-3(n-3)+9n-10+6n  \label{3.7}\\
         &=& f_1(n,4). \notag
\end{eqnarray}
The equality holds in (\ref{3.7}) if and only if $G-u \cong U_{n-1,4}$, $d_G(v)=n-3$, $2\sum_{x \in N(v)\backslash \{u\}} \varepsilon_G(x)+3\sum_{x \in V_G\backslash N[v]}\varepsilon_G(x) = 6n.$ That is to say, $\xi^d(G)=f_1(n,4)$ holds if and only if $G \cong U_{n,4}.$

\medskip
\noindent$\bullet$ $m\ge 5$. We prove the result by induction on $n.$ If $n=2m,$ then the result follows from Theorem 3.1. Suppose that $n > 2m$ and the result holds for graphs in $\mathscr{U}_{k,m}, k\le n-1.$ Let $G \in \mathscr{U}_{n,m}.$ By Lemma 2.4 and the induction hypothesis, it is easily seen that
\begin{eqnarray}
 \xi^d(G)&\geq& \xi^d(G-u)-3d_G(v)+9n-10+2\sum\limits_{x \in N(v))\backslash \{u\}} \varepsilon_G(x)+3\sum\limits_{x \in V_G\backslash N[v]}\varepsilon_G(x) \label{eq:3.8}\\
         &=&  \xi^d(G-u)-9d_G(v)+21n-28 \label{eq:3.9} \\
         &\geq& f_1(n-1,m)-9(n-m+1)+21n-28\ \ \ \ \text{(by induction and (\ref{eq:2.12}))}  \label{eq:3.10}\\
         &=& f_1(n,m). \notag
\end{eqnarray}
Equality in (\ref{eq:3.8}) holds if and only if $\varepsilon_G(v)=2$ and $\varepsilon_{G-u}(x)=\varepsilon_G(x)$ for all $x \in V_{G-u},$ which implies
\begin{eqnarray*}
&& G-u \in \left\{H(n-1,5;1^{k_1+1},1^{k_2+1}S_{t_1,t_2,\dots,t_{m-4}},1^{k_3+1})\right\} \bigcup \left\{H(n-1,4;1^{k_1},1^{k_2}S_{t_1,t_2,\dots,t_{m-3}},1^{k_3})\right\} \\
&& \ \ \ \ \ \ \ \ \ \ \ \ \bigcup\left\{H(n-1,3;1^{k_1+1}S_{t_1,t_2,\dots,t_{m-3}},1^{k_2+1},1^{k_3+1})\right\}\bigcup \left\{H(n-1,5;1^{k_1},1^{k_2}S_{t_1,t_2,\dots,t_{m-3}},1^0)\right\}  \\
&& \ \ \ \ \ \ \ \ \ \ \ \  \bigcup\left\{H(n-1,3;1^{k_1}S_{t_1,t_2,\dots,t_{m-2}},1^{k_2},1^0)\right\},\ \ \text{$k_1,k_2,k_3\geq 0,\, t_i\geq 1, i=1,2,\dots,m-2.$}
\end{eqnarray*}
So we have $\varepsilon_G(x)=4$ for all $x \in V_G\backslash N[v]$ and $\varepsilon_G(x)=3$ for $x \in N(v)$, whence equality in (\ref{eq:3.9}) holds. The equality in (\ref{eq:3.10}) holds if and only if $G-u \cong U_{n-1,m}$, $d_G(v)=n-m+1$. That is to say, $\xi^d(G)= f_1(n,m)$ if and only if $G \cong U_{n,m}$ for $m\ge 5.$

This completes the proof.
\end{proof}
\begin{thm}
Let $G \in \mathscr{U}_{n,2}\backslash\{U_{n,2}\}.$  Then $\xi^d(G) \geq 60,$ for $n=5$, with equality if and only if $G \cong C_5,$ whereas $\xi^d(G) \geq 6n^2-11n-16,$ for $n \geq 6,$ with equality if and only if $G \cong H_{n,4}.$
\end{thm}
\begin{proof}
For the case $n=5,$ our result holds immediately from Table 2. Hence, in what follows we consider $n\geq 6.$ Combining (i) of Theorem 3.3 yields that $\mathscr{U}_{n,2}\backslash\{U_{n,2}\}=\{C_3(1^{k_1},1^{k_2},1^0):k_1+k_2=n-3,k_1>k_2\geq 1\}\cup \{H_{n,4}\}\cup \{H(n,3;1^0S_{n-4},1^0,1^0)\}.$

From \cite{12} we know that  $C_3(1^{n-4},1^1,1^0)$ is the only graph with the second minimal EDS among the $n$-vertex unicyclic graphs of girth $3$. And it is easy to check that
\begin{eqnarray*}
 && \xi^d(C_3(1^{n-4},1^{1},1^0))-\xi^d(H_{n,4})= (6n^2-11n-15)-(6n^2-11n-16)=-1<0.
\end{eqnarray*}
Then $H_{n,4}$ is the graph with the second minimal EDS among $\mathscr{U}_{n,2},$ as desired.
\end{proof}

\begin{thm}
Let $G \in \mathscr{U}_{n,3}$ with $n \geq 6.$
\begin{wst}
\item[{\rm (i)}]If $n=6,7,$ then $U_{n,3}^*$ is the unique unicyclic graph with the second minimal EDS among the graphs in $\mathscr{U}_{n,3}.$
\item[{\rm (ii)}]If $n=8,9,$ then $U_{n,3}$ is the unique unicyclic graph with the second minimal EDS among the graphs in $\mathscr{U}_{n,3}.$
\item[{\rm (iii)}]If $n \geq 10$, then $\xi^d(G) \geq 6n^2-5n-53$ for $G \in \mathscr{U}_{n,3}\backslash \{U_{n,3}\},$ with equality if and only if $G \cong U_{n,3}'$.
\end{wst}
\end{thm}
\begin{proof}
(i) and (ii) follow directly from Table 2.

(iii)\ \ For convenience, let $g_3(n)=6n^2-5n-53$ with $n\ge 10$. It is routine to check that $G \ncong C_n$. We prove our result by induction on $n.$

If $n=10$, then our result follows immediately from Table 2. So we consider that $n \geq 11$ and assume that our result holds for graphs in $\mathscr{U}_{k,3}\backslash \{U_{k,3}\}$ with $k \leq n-1.$  Hence, for graph $G\in \mathscr{U}_{n,3}$, in view of Lemma 2.4, we have inequality holds in (\ref{eq:2.4}), with equality if and only if $\varepsilon_G(v)=2$ and $\varepsilon_{G-u}(x)=\varepsilon_G(x)$ for all $x \in V_{G-u},$ where $u$ is a pendant vertex and $v$ is the unique neighbor of $u$.
Hence, by the definition of $\mathscr{U}^2_{n,3},$ we have
\begin{eqnarray*}
&&G\in \mathscr{U}^2_{n,3}= \left\{H(n,5;1^{k_1},1^{k_2+1}2^0,1^0)\right\}\bigcup \left\{H(n,4;1^{k_1+1},1^{k_2+1}2^0,1^{k_3})\right\}\bigcup \left\{H(n,4;1^{k_1},1^{k_2}S_{t_1},1^{k_3})\right\}\\
                 & & \ \ \ \ \ \ \ \ \ \ \ \ \  \,\,\,\,\,\,\,\, \bigcup\left\{H(n,3;1^{k_1+1}2^0,1^{k_2+1},1^{k_3+1})\right\}\bigcup \left\{H(n,3;1^{k_1}S_{t_1},1^{k_2},1^0)\right\},\   k_1, k_2, k_3\geq 0, t_1 \geq 1.
\end{eqnarray*}
In view of (\ref{eq:2.4}), let $\phi(G)=2\sum_{x \in N(v)\backslash \{u\}} \varepsilon_G(x)+3\sum_{x \in V_G\backslash N[v]}\varepsilon_G(x).$ By the structure of graphs in $\mathscr{U}^2_{n,3}\setminus \{U_{n,3}'\}$ and by direct calculation, we have
\begin{align*}
\phi(H(n,5;1^{k_1},1^{k_2+1}2^0,1^0))\geq& \ \phi(H(n,5;1^0,1^{n-5},1^0))=\phi(U_{n,3}')=6n-10,\\
\phi(H(n,4;1^{k_1+1},1^{k_2+1}2^0,1^{k_3})) \geq&\ \phi(H(n,4;1^1,1^{n-5},1^0))=6n-8 ,\\
\phi(H(n,4;1^{k_1},1^{k_2}S_{t_1},1^{k_3}))\geq&\ \phi(H(n,4;1^0,1^{n-6}2^1,1^0))=6n,\\
\phi(H(n,3;1^{k_1+1}2^0,1^{k_2+1},1^{k_3+1}))\geq &\ \phi(H(n,3;1^{n-5}2^0,1^1,1^1))=6n-10,\\
\phi(H(n,3;1^{k_1}S_{t_1},1^{k_2},1^0))\geq& \ \phi(H(n,3;1^{n-6}2^1,1^1,1^0))=6n.
\end{align*}
Hence, in view of (\ref{eq:2.4}) and by induction we have
\begin{eqnarray}
\text{RHS of (\ref{eq:2.4})} &\geq& g_3(n-1)-3(n-3)+9n-10+6n-10 \ \ \ \ \text{(by induction and (\ref{eq:2.13}))} \label{3.12}\\
         &=& g_3(n). \notag
\end{eqnarray}
The equality holds in (\ref{3.12}) if and only if $G-u \cong U_{n-1,3}'$, $d_G(v)=n-3$ and $ 2\sum_{x \in N(v))\backslash \{u\}} \varepsilon_G(x)\linebreak +3\sum_{x \in V_G\backslash N[v]}\varepsilon_G(x) = 6n -10.$ That is to say, $\xi^d(G)=g_3(n)$ if and only if $G \cong U_{n,3}'.$

This completes the proof.
\end{proof}

\begin{thm}
Let $G \in \mathscr{U}_{n,4}\backslash \{U_{n,4}\},$ where $n \geq 9.$
\begin{wst}
\item[{\rm (i)}]If $9 \leq n  \leq15$ then $\xi^d(G) \geq 6n^2+17n-147$ with equality if and only if $G \cong C_5(1^1,1^{n-7},1^1).$
\item[{\rm (ii)}]If $n \geq 16$, then $\xi^d(G) \geq 6n^2+14n-100$ with equality if and only if $G \cong U_{n,4}'.$
\end{wst}
\end{thm}
\begin{proof}
(i) follows immediately from Table 2.

(ii)\ \ Suppose $n \geq 16$ and let $g_4(n)=6n^2+14n-100$. It is routine to check that $G \ncong C_n$. We prove our result by induction on $n.$

If $n=16$, then our result follows immediately from Table 2. So we consider that $n \geq 17$ and assume that our result holds for graphs in $\mathscr{U}_{k,4}\backslash \{U_{k,4}\}$ with $k \leq n-1.$  Hence, for graph $G\in \mathscr{U}_{n,4}$, in view of Lemma 2.4, we have inequality holds in (\ref{eq:2.4}), with equality if and only if $\varepsilon_G(v)=2$ and $\varepsilon_{G-u}(x)=\varepsilon_G(x)$ for all $x \in V_{G-u},$ where $u$ is a pendant vertex and $v$ is the unique neighbor of $u$.
Hence, by the definition of $\mathscr{U}^2_{n,4},$ we obtain that $G$ is in
\begin{eqnarray*}
&&\mathscr{U}^2_{n,4}= \left\{H(n,5;1^{k_1+1},1^{k_2+1}2^0,1^{k_3+1})\right\}\bigcup \left\{H(n,5;1^{k_1},1^{k_2}S_{t_1},1^{k_3})\right\}\bigcup \left\{H(n,4;1^{k_1+1},1^{k_2+1}S_{t_1},1^{k_3})\right\}\\
                 & & \ \ \ \ \ \ \ \ \ \ \bigcup\left\{H(n,4;1^{k_1},1^{k_2}S_{t_1,t_2},1^{k_3})\right\}\bigcup \left\{H(n,3;1^{k_1+1}S_{t_1},1^{k_2+1},1^{k_3+1})\right\}\bigcup \{H(n,3;1^{k_1+1}S_{t_1,t_2},1^{k_2},1^{k_3})\},
\end{eqnarray*}
where $k_1,k_2,k_3 \geq 0$ and $t_1,t_2\geq 1.$  In view of (\ref{eq:2.4}), let $\phi(G)=2\sum_{x \in N(v)\backslash \{u\}} \varepsilon_G(x)+3\sum_{x \in V_G\backslash N[v]}\varepsilon_G(x).$ By the structure of graphs in $\mathscr{U}^2_{n,4}\setminus \{U_{n,4}\}$ and by direct calculation, we have
\begin{align*}
\phi(H(n,5;1^{k_1+1},1^{k_2+1}2^0,1^{k_3+1}))\geq& \ \phi(H(n,5;1^1,1^{n-7}2^0,1^1))=6n+6,\\
\phi(H(n,5;1^{k_1},1^{k_2}S_{t_1},1^{k_3})) \geq&\ \phi(H(n,5;1^0,1^{n-7}2^1,1^0))=6n+6,\\
\phi(H(n,4;1^{k_1+1},1^{k_2+1}S_{t_1},1^{k_3}))\geq& \ \phi(H(n,4;1^1,1^{n-7}2^1,1^0))=6n+6,\\
\phi(H(n,4;1^{k_1},1^{k_2}S_{t_1,t_2},1^{k_3}))\geq& \ \phi(H(n,4;1^0,1^{n-8}2^2,1^0))=6n+6,\\
\phi(H(n,3;1^{k_1+1}S_{t_1},1^{k_2+1},1^{k_3+1}))\geq& \ \phi(H(n,3;1^{n-7}2^1,1^1,1^1))=6n+6,\\
\phi(H(n,3;1^{k_1+1}S_{t_1,t_2},1^{k_2},1^{k_3}))\geq& \ \phi(H(n,3;1^{n-8}2^2,1^1,1^0))=6n+6.
\end{align*}
Hence, in view of (\ref{eq:2.4}) and by induction we have
\begin{eqnarray}
\text{RHS of (\ref{eq:2.4})} &\geq& g_4(n-1)-3(n-4)+9n-10+6n+6 \ \ \ \ \text{(by induction and (\ref{eq:2.13}))} \label{3.13}\\
         &=& g_4(n). \notag
\end{eqnarray}
The equality holds in (\ref{3.13}) if and only if $G-u \cong U_{n-1,4}'$, $d_G(v)=n-4$ and $ 2\sum_{x \in N(v))\backslash \{u\}} \varepsilon_G(x)\linebreak +3\sum_{x \in V_G\backslash N[v]}\varepsilon_G(x) = 6n+6.$ That is to say, $\xi^d(G)=g_4(n)$ if and only if $G \cong U_{n,4}'.$

This completes the proof.
\end{proof}

\begin{thm}
Let $G \in \mathscr{U}_{n,m} \backslash \{U_{n,m}\}$ with $m \geq 5$, then $\xi^d(G) \geq 6n^2 + m^2 +9mn - 28m - 22n -4$ with equality if and only if $G \cong U_{n,m}'$.
\end{thm}
\begin{proof}
For convenience, let $f_2(n,m)=6n^2 + m^2 +9mn - 28m - 22n -4$. We prove the result by induction on $n.$

If $n=2m,$ then our result follows from Theorem 3.2. Suppose that $n > 2m$ and the result holds for graphs in $\mathscr{U}_{k-1,m}\backslash \{U_{k-1,m}\}$ for
$k\le n-1$. Now we consider $G \in \mathscr{U}_{n,m}\backslash \{U_{n,m}\}.$ By Lemma 2.4 and the induction hypothesis, there exists a pendant vertex $u$ with a unique neighbor $v$ such that
\begin{eqnarray}
 \xi^d(G)&\geq& \xi^d(G-u)-3d_G(v)+9n-10+2\sum\limits_{x \in N(v))\backslash \{u\}} \varepsilon_G(x)+3\sum\limits_{x \in V_G\backslash N[v]}\varepsilon_G(x) \label{3.14}\\
         &=&  \xi^d(G-u)-9d_G(v)+21n-28 \label{3.15} \\
         &\geq& f_2(n-1,m)-9(n-m)+21n-28 \ \ \ \ \text{(by induction and (\ref{eq:2.13}))} \label{3.16}\\
         &=& f_2(n,m). \notag
\end{eqnarray}
Equality in (\ref{3.14}) holds if and only if $\varepsilon_G(v)=2$ and $\varepsilon_{G-u}(x)=\varepsilon_G(x)$ for all $x \in V_{G-u},$ which implies that
\begin{eqnarray*}
&& G-u \in \{H(n-1,5;1^{k_1+1},1^{k_2+1}S_{t_1,t_2,\dots,t_{m-4}},1^{k_3+1})\} \bigcup \{H(n-1,4;1^{k_1},1^{k_2}S_{t_1,t_2,\dots,t_{m-3}},1^{k_3})\} \\
&& \ \ \ \ \ \ \ \ \ \ \ \ \bigcup\{H(n-1,3;1^{k_1+1}S_{t_1,t_2,\dots,t_{m-3}},1^{k_2+1},1^{k_3+1})\}\bigcup \{H(n-1,5;1^{k_1},1^{k_2}S_{t_1,t_2,\dots,t_{m-3}},1^0)\} \\
&& \ \ \ \ \ \ \ \ \ \ \ \ \bigcup \{H(n-1,3;1^{k_1}S_{t_1,t_2,\dots,t_{m-2}},1^{k_2},1^0)\},\ \ k_1,k_2,k_3\geq 0,\, t_i\geq 1, i=1,2,\dots,m-2.
\end{eqnarray*}
So we have $\varepsilon_G(x)=4$ for all $x \in V_G\backslash N[v]$ and $\varepsilon_G(x)=3$ for $x \in N(v)$, whence equality in (\ref{3.15}) holds.
The equality in (\ref{3.16}) holds if and only if $G-u \cong U_{n-1,m}'$, $d_G(v)=n-m$. That is to say, $\xi^d(G)=f_2(n,m)$ if and only if $G \cong U_{n,m}'.$

This completes the proof.
\end{proof}


\begin{thebibliography}{99}
\small \setlength{\itemsep}{.2mm}

\bibitem{0}A.R. Ashrafi, M. Saheli, M. Ghorbani, The eccentric connectivity index of nanotubes and nanotori,
J. Comput. Appl. Math. 235 (2011) 4561-4566.

\bibitem{D-I}J.A. Bondy, U.S.R. Murty, Graph Theory with Applications. Macmillan Press, New York (1976)

\bibitem{C-R-S1}A, Chang. F, Tian, On the spectral radius of unicyclic graphs with perfect matching, Linear Algebra
Appl. 370 (2003) 237-250.



\bibitem{4}A. Dobrynin, R. Entringer, I. Gutman, Wiener index of trees: theory and applications, Acta Appl. Math. 66 (2001) 211-249.

\bibitem{05}A. Dobrynin, A.A. Kochetova, Degree distance of a graph: A degree analogue of the Wiener index, J. Chem. Inf. Comput. Sci. 34 (1994) 1082-1086.



\bibitem{S0}S. Gupta, M. Singh, A.K. Madan, Eccentric distance sum: A novel graph invariant for predicting biological and physical properties, J. Math. Anal.
Appl. 275 (2002) 386-401.



\bibitem{1}I. Gutman, O.E. Polansky, Mathematical Concepts in Organic Chemistry, Springer, Berlin (1986)

\bibitem{09}I. Gutman, Selected properties of the Schultz molecular topological index, J. Chem. Inf. Comput. Sci. 34 (1994) 1087-1089.

\bibitem{2}Y.P. Hou, J.S. Li, Bounds on the largest eigenvalues of trees with a given size of matching, Linear
Algebra Appl. 342 (2002) 203-217.

\bibitem{H-H-X-K}H.B. Hua, K.X. Xu, W.N. Shu, A short and unified proof of Yu et al.'s two results on the eccentric
distance sum, J. Math. Anal. Appl. 382 (2011) 364-366.

\bibitem{3}H.B. Hua, S.G. Zhang, K.X. Xu, Further rasults on the eccentric distance
sum,  Discrete Appl. Math. 160 (2012) 170-180.




\bibitem{13}A. Ili\'c, G.H. Yu, L.H. Feng, On the eccentric distance sum of graphs, J. Math. Anal. Appl. 381 (2011) 590-600.

\bibitem{17}S.C. Li, M. Zhang, G.H. Yu, L.H. Feng, On the extremal values of the eccentric distance sum of trees, J. Math. Anal. Appl. 390 (2012) 99-112.



\bibitem{19}B. McKay, Nauty \url{http://cs.anu.edu.au/~bdm/nauty/}




\bibitem{12}G.H. Yu, L.H. Feng, A. Ili\'c, On the eccentric distance sum of trees and unicyclic graphs,
J. Math. Anal. Appl. 375 (2011) 99-107.
\end{thebibliography}
\end{document}